\documentclass[12pt,draftcls,onecolumn]{IEEEtran}
\usepackage{amsfonts}
\usepackage{amssymb}
\usepackage{amsmath}
\usepackage{mathptmx}
\usepackage{color}
\usepackage{cancel}
\usepackage{bbm}
\usepackage{psfrag}
\usepackage{mathrsfs}
\def\be{\begin{equation}}
\def\ee{\end{equation}}
\def\bea{\begin{eqnarray}}
\def\eea{\end{eqnarray}}
\def\beann{\begin{eqnarray*}}
\def\eeann{\end{eqnarray*}}

\def\ns{\hspace{-1mm}}
\newcommand{\real}{{\mathbb{R}}}

\def\spacingset#1{\def\baselinestretch{#1}\small\normalsize}
\spacingset{1.3}

\newtheorem{lemma}{Lemma}
\newtheorem{theorem}{Theorem}
\newtheorem{remark}{Remark}
\newtheorem{corollary}{Corollary}
\newtheorem{proposition}{Proposition}

\newtheorem{problem}{Problem}

\def\be{\begin{equation}}
\def\ee{\end{equation}}
\def\bea{\begin{eqnarray}}
\def\eea{\end{eqnarray}}
\def\beann{\begin{eqnarray*}}
\def\eeann{\end{eqnarray*}}

\def\ns{\hspace{-1mm}}

\def\proof{\noindent{\bf{\em Proof:}\ \ }}
\def\QED{\mbox{\rule[0pt]{1.5ex}{1.5ex}}}
\def\endproof{\hspace*{\fill}~\QED\par\endtrivlist\unskip}

\newcommand{\ima}{\operatorname{im}}

\newcommand{\defi}{\stackrel{\text{\tiny def}}{=}}

\def\tra{\top}

\def\gR{{\cal R}}

\def\gV{{\cal V}}

\def\gS{{\cal S}}

\def\bmat{\left[ \begin{array}}
\def\emat{\end{array} \right]}

\def\bmat{\left[ \begin{array}}
\def\emat{\end{array} \right]} 
\def\bsmat{\left[ \begin{smallmatrix}}
\def\esmat{\end{smallmatrix} \right]} 
\newcommand{\spanR}{\operatorname{span}}

\def\gR{{\cal R}}

\def\gV{{\cal V}}

\def\gS{{\cal S}}

\begin{document}
\title{\LARGE{{Continuous-Time Singular Linear-Quadratic Control: Necessary and Sufficient Conditions for the Existence of Regular Solutions}}}

\author{Augusto~Ferrante,
        Lorenzo~Ntogramatzidis
\thanks{Augusto Ferrante is with the Dipartimento di Ingegneria dell'Informazione,
         Universit\`a di Padova,
         via Gradenigo, 6/B -- I-35131 Padova, Italy.
        E-mail: {\tt augusto@dei.unipd.it.} Research carried out while visiting Curtin University, Australia.}
\thanks{L. Ntogramatzidis is with the Department of Mathematics and
Statistics, Curtin University, Perth WA 6845,
Australia. E-mail: {\tt L.Ntogramatzidis@curtin.edu.au. }}

\thanks{Partially supported by the Italian Ministry for Education and Research (MIUR) under PRIN grant n. 20085FFJ2Z.}

}
 
\markboth{DRAFT}{Shell \MakeLowercase{\textit{et al.}}: Bare Demo of IEEEtran.cls for Journals}

\maketitle

\vspace{-1cm}

\begin{abstract}
The purpose of this paper is to close the remaining gaps in the understanding of the role that the constrained generalized continuous algebraic Riccati equation plays in singular linear-quadratic (LQ) optimal control. 
{Indeed, in spite of the vast literature on  LQ problems, it is only in a recent paper that a sufficient condition for the existence of a non-impulsive optimal control has for the first time connected this equation with the singular LQ optimal control problem.
In this paper, we establish four equivalent conditions providing a complete picture that connects the singular LQ problem  with the generalized continuous algebraic Riccati equation and with the geometric 
properties of the underlying system.}
\end{abstract}

\IEEEpeerreviewmaketitle

\section{Introduction}
\label{secintro} 
This paper addresses the continuous-time linear quadratic (LQ) optimal control problem when the matrix weighting the input in the cost function, traditionally denoted by $R$, is possibly singular. 
 This problem has a long history. It has been investigated in several papers and with the use of different techniques, see
\cite{Hautus-S-83,Willems-KS-86,Saberi-S-87,Prattichizzo-MN-04,Kalaimani-BC-13} and the references cited therein.
 In particular, in the classical contributions
\cite{Hautus-S-83} and \cite{Willems-KS-86} it was proved that {\bf i)} an optimal solution of the singular LQ problem exists for all initial conditions if the class of allowable controls is extended to include
distributions;  {\bf ii)} the regular part of the optimal control can still be written as a static state feedback $u=-K\,x$ as in the regular case.
In the discrete time,  the solution of regular and singular finite and infinite-horizon LQ problems can be found resorting to the so-called  {\em constrained generalized discrete algebraic Riccati equation}, see  \cite{Ferrante-N-13-1,Ferrante-N-12-sub} and also \cite{Stoorvogel-S-98}.  A similar generalization has been carried out for the continuous-time algebraic Riccati equation in \cite{Ionescu-O-96-1}, where the {\em 
constrained generalized Riccati equation} was defined in such a way that the inverse of $R$ appearing in the standard Riccati equation is replaced by its pseudo-inverse. 
On the other hand, until very recently this counterpart of the generalized discrete algebraic Riccati equation was only studied without any understanding of its links with the linear quadratic optimal control problem. 

The recent paper \cite{Ferrante-N-14} was the first attempt to provide a description of the role played by the constrained generalized continuous algebraic Riccati equation  in singular LQ optimal control problems. Such role does not trivially follow from the analogy with the discrete case, as one can immediately realize by considering the fact that in the continuous time, whenever the optimal control involves distributions, none of the solutions of the constrained generalized Riccati equation is optimizing. In particular, in \cite{Ferrante-N-14}  it was shown that when the continuous-time constrained generalized Riccati equation possesses a symmetric solution, the corresponding LQ problem admits {a {\em regular} (i.e. impulse-free)} solution, and an optimal control can always be expressed as a state-feedback.
{
This is just a single trait of a rich picture where necessary and sufficient conditions for the existence of regular solutions are given in terms of the algebraic and  geometric structures of the underlying system. 
In particular, the algebraic   structure refers to the existence of solutions to the associated generalized  algebraic Riccati equation. The purpose of this paper is to provide a full illustration of this picture.  
}\\[-5mm]

{\bf Notation.}  The image and the kernel of matrix $M$ are denoted by $\ima\,M$ and $\ker\,M$, respectively, while the transpose and the Moore-Penrose pseudo-inverse of $M$ are denoted by $M^\tra$ and $M^\dagger$, respectively. Given a quadruple of matrices $(A,B,C,D)$, where $A\in \real^{n \times n}$, $B\in \real^{n \times m}$, 
$C\in \real^{p \times n}$ and $D\in \real^{p \times m}$, {we denote by $\gV^\star$ the largest output-nulling subspace,  by $\gS^\star$ the smallest input containing subspace, and  by $\gR^\star$ the largest reachability output-nulling subspace, see \cite{Trentelman-SH-01} for details.} 

\subsection{Preliminaries}
{A key role in this paper will be played by} the following matrix equation
\begin{equation}
 X\,A+A^\tra\,X-(S+X\,B)\,R^{\dagger}\,(S^\tra\!+B^\tra X)+Q=0, \label{gcare}
\end{equation}
with $Q, A\in \real^{n \times n}$, $B,S \in \real^{n \times m}$, $R \in \real^{m \times m}$ and we make the following standing assumption:
\begin{equation}
\label{sd}
\Pi \defi \bmat{cc}  Q & S \\ S^\tra & R \emat=\Pi^\tra \ge 0.
\end{equation}
Thus, the {\em Popov matrix} $\Pi$ can be factorized in terms of two matrices $C\in \real^{p \times n}$ and $D\in \real^{p \times m}$ as
\bea
\label{fact}
\Pi=\bmat{cc} C^\top \\ D^\top \emat \bmat{cc} C & D \emat.
\eea

Let us identify $\Sigma$ with the triple $(A,B,\Pi)$. 
Eq. (\ref{gcare}) is often referred to as the {\em generalized continuous algebraic Riccati equation} GCARE($\Sigma$), and represents a generalization of the classic continuous algebraic Riccati equation CARE($\Sigma$)
\begin{equation}
 X\,A+A^\tra\,X-(S+X\,B)\,R^{-1}\,(S^\tra\!+B^\tra X)+Q=0, \label{care}
\end{equation}
arising in infinite-horizon LQ problems since in the present setting $R$ is allowed to be singular. Eq. (\ref{gcare}) along with the additional condition
\begin{equation}
 \ker R \subseteq \ker (S+X\,B), \label{kercond}
\end{equation}
is usually referred to as {\em constrained generalized continuous algebraic Riccati equation}, and is denoted by CGCARE($\Sigma$). Observe that from (\ref{sd}) we have $\ker R \subseteq \ker S$, which implies that (\ref{kercond}) is equivalent to $\ker R \subseteq \ker (X\,B)$. 
 
{
The classic LQ optimal control problem can be stated as follows
 \begin{problem}\label{prb1} Find a control input $u(t)$, $t \ge 0$, that minimizes the performance index 
   \begin{eqnarray}
 \label{costinf}
J_\infty(x_{\scriptscriptstyle 0},u)\ns&\ns=\ns&\ns\int_{\scriptscriptstyle 0}^\infty \bmat{cc} x^\tra(t) & u^\tra(t) \emat \bmat{cc} Q & S \\ S^\tra & R \emat \bmat{c} x(t) \\ u(t) \emat\,dt
\end{eqnarray}
subject to the constraint 
 \begin{equation}
 \label{eqsys}
 \dot{x}(t)=A\,x(t)+B\,u(t), \qquad x(0)=x_{\scriptscriptstyle 0} \in \real^n.
 \end{equation}
\end{problem}
} 
{
We consider $u$ to be a solution of Problem \ref{prb1} only if the corresponding value of the performance index is finite.}\footnote{We make this remark since, if the cost is unbounded for every control, one {might alternatively say} that all controls are optimal since they all lead to the same value of the performance index.}
{Moreover, we say that a solution $u^\ast$ of Problem \ref{prb1} is regular if $u^\ast\in{\mathcal C}_\infty[0,\infty)$.}
 
 It is well-known that when $R$ is positive definite, the optimal control (when it exists) does not include distributions, since in such a case an impulsive control $u$ will always cause $J_\infty(x_{\scriptscriptstyle 0},u)$ to be unbounded for any $x_{\scriptscriptstyle 0}\in \real^n$. If $R$ is only positive semidefinite, in general the optimal solution can contain distributions, given by Dirac delta distributions and its derivatives.

\section{Main result}

The main result of this paper is the following theorem, whose proof will be developed in several steps in the sequel.
 
 \begin{theorem}
 \label{main}
 The following statements are equivalent:
 \begin{description}
 \item{{\bf \em \hspace{-1cm}  (A).}} {For every initial state $x_{\scriptscriptstyle 0}\in \real^n$, Problem \ref{prb1} admits a regular solution;}
 \item{{\bf \em \hspace{-1cm}  (B).}}  There exists a symmetric and positive semidefinite solution of CGCARE($\Sigma$);
 \item{{\bf \em \hspace{-1cm}  (C).}} There exists a symmetric solution of CGCARE($\Sigma$), and for {each initial  state} $x_{\scriptscriptstyle 0}\in \real^n$, there exists $u_{\scriptscriptstyle 0}(t)$ such that $J_\infty(x_{\scriptscriptstyle 0},u_{\scriptscriptstyle 0})$ is finite;
  \item{{\bf \em \hspace{-1cm}  (D).}} For any factorization (\ref{fact}), the subspaces $\gS^\star$ and $\gR^\star$ of the quadruple $(A,B,C,D)$ coincide, and and for {each initial state} $x_{\scriptscriptstyle 0}\in \real^n$, there exists $u_{\scriptscriptstyle 0}(t)$ such that $J_\infty(x_{\scriptscriptstyle 0},u_{\scriptscriptstyle 0})$ is finite.
 \end{description}
 \end{theorem}

 \begin{remark}{
 {Existence, for each $x_{\scriptscriptstyle 0}$, of a control function  $u_{\scriptscriptstyle 0}(t)$ such that $J_\infty(x_{\scriptscriptstyle 0},u_{\scriptscriptstyle 0})$ is finite, is a very natural and mild  condition. Its testability, however, is not obvious.
 It has been shown in \cite{Geerts-H-89} that such condition is equivalent to the following neat and easily testable geometric condition:} 
 \[
 \gV^\star+\langle A,\ima B \rangle+{\cal X}_{\rm stab}=\real^n,
 \]
 where $\gV^\star$ is the largest output-nulling subspace of the quadruple $(A,B,C,D)$, $\langle A,\ima B \rangle$ is the reachable subspace (i.e., the smallest $A$-invariant subspace containing the range of $B$), {and ${\cal X}_{\rm stab}$ is the  $A$-invariant subspace  corresponding to the asymptotically stable uncontrollable eigenvalues of $A$} (so that, in other words, the sum $\langle A,\ima B \rangle+{\cal X}_{\rm stab}$ is the stabilizable subspace of the pair $(A,B)$).
 }
 \end{remark}

\section{Ancillary results and proof of main result}
 
 The following notation is used throughout the paper. 
We denote by $G \defi I_m-R^\dagger R$   the orthogonal projector that projects onto $\ker R$.  Moreover, we consider a non-singular matrix $T=[T_1\mid T_2]$ where
$\ima T_1=\ima R$ and $\ima T_2=\ima G$, and we define $B_1\defi BT_1$ and $B_2 \defi BT_2$. Finally, to any $X=X^\tra \in \real^{n \times n}$  we associate\begin{eqnarray}
Q_X& \defi & Q+A^\tra X+X\,A, \qquad
S_X   \defi   S+X\, B, \label{defgx} \\
K_X & \defi & R^\dagger\, (S^\tra+B^\tra\,X)=R^\dagger\, S_X^\tra, \qquad A_X \defi  A-B\,K_X, \\
  \Pi_X  & \defi & \left[ \begin{array}{cc} Q_X & S_X \\ S_X^\tra & R \end{array} \right]. \label{KX}
 \end{eqnarray}

 The following {result}, which is the main result of \cite{Ferrante-N-14}, establishes that when CGCARE($\Sigma$) admits at least one symmetric solution, and the performance index can be rendered finite with a certain control function for every initial state, the corresponding LQ optimal control problem admits impulse-free controls.

 \begin{proposition}
 \label{theprev}
 Suppose CGCARE($\Sigma$) admits symmetric solutions, and that for every $x_{\scriptscriptstyle 0}$ there exists an input $u(t) \in \real^m$, with $t\ge0$, such that
$J_\infty(x_{\scriptscriptstyle 0},u)$ in (\ref{costinf}) is finite. Then:
\begin{itemize} 
 \item A solution $\overline{X}=\overline{X}^\tra\ge 0$ of CGCARE($\Sigma$) is obtained as the limit of the time varying matrix generated by integrating 
 \begin{eqnarray}
\dot{X}(t) =X(t)\,A+A^\tra\,X(t) -(S+X(t)\,B)\,R^{\dagger}\,(S^\tra+B^\tra X(t))+Q, \label{grde111}
\end{eqnarray}

  with the zero initial condition $X(0)=0$.
 \item The value of the optimal cost is {\em $x_{\scriptscriptstyle 0}^\tra \,\overline{X} \,x_{\scriptscriptstyle 0}$}.
 \item $\overline{X}$ is the minimum positive semidefinite solution  of CGCARE($\Sigma$).
 \item The set of {\em all} optimal controls minimizing {the cost} in (\ref{costinf}) can be parameterized as 
 \be
u(t)=-R^\dagger S_{\scriptscriptstyle \overline{X}}^\tra\,x(t)+G\,v(t), \label{optcontrinf}
\ee
with arbitrary $v(t)$.
\end{itemize}
\end{proposition}

It is easy to see that Proposition \ref{theprev} proves that the implications {\bf {\em (C)}} $\Rightarrow$ {\bf {\em (B)}} and {\bf {\em (C)}} $\Rightarrow$ {\bf {\em (A)}} in Theorem \ref{main} hold true.
The following Proposition shows that {\bf {\em (B)}} $\Rightarrow$ {\bf {\em (C)}} as well.


\begin{proposition}
If there exists a symmetric positive semidefinite solution $\overline{X}=\overline{X}^\top \ge0$ of CGCARE($\Sigma$), then for all initial states $x_{\scriptscriptstyle 0}\in \real^n$, there exists $u_{\scriptscriptstyle 0}(t)$ such that $J_\infty(x_{\scriptscriptstyle 0},u_{\scriptscriptstyle 0})$ is finite.
\end{proposition}
\proof
Let $u_{\scriptscriptstyle 0}(t) = -R^\dagger\,S_{\scriptscriptstyle \overline{X}}^\top\,x(t)$, where we recall that $S_{\scriptscriptstyle \overline{X}}=S+{\overline{X}}\,B$.
We can write the state equation as
\[
\dot{x}(t)=A_{\scriptscriptstyle \overline{X}}\,x(t), 
\]
where $A_{\scriptscriptstyle \overline{X}}=A-B\,R^\dagger\,S_{\scriptscriptstyle \overline{X}}^\top$. This obviously implies that 
$x(t)=e^{A_{\scriptscriptstyle \overline{X}}\,t}\,x_{\scriptscriptstyle 0}$. We have
\beann
J_\infty(x_{\scriptscriptstyle 0},u_{\scriptscriptstyle 0}) \ns&\ns = \ns&\ns \int_{\scriptscriptstyle 0}^\infty x^\top(t) \bmat{cc} I_n & -S_{\scriptscriptstyle \overline{X}}\,R^\dagger \emat 
\bmat{cc}  Q & S \\ S^\tra & R \emat \bmat{c} I_n \\ -R^\dagger \,S_{\scriptscriptstyle \overline{X}}^\top \emat x(t) \,dt \\
\ns&\ns = \ns&\ns \int_{\scriptscriptstyle 0}^\infty x^\top(t) \left[ Q-S_{\scriptscriptstyle \overline{X}}\,R^\dagger S_{\scriptscriptstyle \overline{X}}^\top+S_{\scriptscriptstyle \overline{X}}\,R^\dagger\,B^\top\,{\overline{X}}+{\overline{X}}\,B\,R^\dagger\,S_{\scriptscriptstyle \overline{X}}^\top \right]\,x(t)\,dt \\
\ns&\ns = \ns&\ns \int_{\scriptscriptstyle 0}^\infty x^\top(t) \left[ -{\overline{X}}\,A-A^\top\,\overline{X}+S_{\scriptscriptstyle \overline{X}}\,R^\dagger\,B^\top\,{\overline{X}}+{\overline{X}}\,B\,R^\dagger\,S_{\scriptscriptstyle \overline{X}}^\top  \right]\,x(t)\,dt \\
\ns&\ns = \ns&\ns -\int_{\scriptscriptstyle 0}^\infty x_{\scriptscriptstyle 0}^\top\,e^{A_{\scriptscriptstyle \overline{X}}^\top\,t} \left[ -{\overline{X}}\,A_{\scriptscriptstyle \overline{X}}-A_{\scriptscriptstyle \overline{X}}^\top\,{\overline{X}}  \right]\,e^{A_{\scriptscriptstyle \overline{X}}\,t} x_{\scriptscriptstyle 0}\,dt \\
\ns&\ns = \ns&\ns \lim_{T \to \infty} \int_{\scriptscriptstyle 0}^T x_{\scriptscriptstyle 0}^\top\,\frac{d}{dt} \left[ -e^{A_{\scriptscriptstyle \overline{X}}^\top\,t} {\overline{X}}\,e^{A_{\scriptscriptstyle \overline{X}}\,t}\right]\, x_{\scriptscriptstyle 0}\,dt \\
\ns&\ns = \ns&\ns \lim_{T \to \infty} x_{\scriptscriptstyle 0}^\top\, \left[ {\overline{X}}-e^{A_{\scriptscriptstyle \overline{X}}^\top\,T} {\overline{X}}\,e^{A_{\scriptscriptstyle \overline{X}}\,T}\right]\,x_{\scriptscriptstyle 0} \le x_{\scriptscriptstyle 0}^\top\,{\overline{X}}\,x_{\scriptscriptstyle 0}.
\eeann
\endproof

The classical papers on singular LQ optimal control \cite{Hautus-S-83,Willems-KS-86} make the strong assumption of stabilizability of the pair $(A,B)$, even when the problem is formulated without a stability constraint on the state trajectory, just to the end of ensuring the convergence of the integral in the cost function. We want to remove this conservative assumption, and only ask for the very weak requirement that there exists a control function that renders the value of the cost function finite. The following {classical result accomplishes this task (we include, for the sake of completeness a very direct proof of this  result).}

\begin{lemma}
\label{lemW}
Consider a regular LQ problem, i.e., with $R=R^\top >0$. If for every $x_{\scriptscriptstyle 0}\in \real^n$ there exists a control function $u(t)\in \real^m$, with $t \ge 0$, such that $J_{\infty}(x_{\scriptscriptstyle 0},u)$ is finite, {then there exist solutions $X=X^\top\geq 0$ of CARE($\Sigma$).
Among such solutions there is a minimal one $\overline{X}$ and the optimal control is given by $u^\ast(t)=-R^{-1}(S^\top+B^\top\,\overline{X})\,x(t)$.}
\end{lemma}
\proof
{Consider the finite-horizon performance index
 \begin{eqnarray}
 J_{T}(x_{\scriptscriptstyle 0},u)&=&\int_{\scriptscriptstyle 0}^T \bmat{cc} x^\tra(t) & u^\tra(t) \emat \bmat{cc} Q & S \\ S^\tra & R \emat \bmat{c} x(t) \\ u(t) \emat\,dt,  \label{cost1}
 \end{eqnarray}
and the Riccati differential equation
\begin{eqnarray}
 \dot{P}_{\scriptscriptstyle T}(t)+P_{\scriptscriptstyle T}(t)\,A+A^\tra\,P_{\scriptscriptstyle T}(t)-(S+P_{\scriptscriptstyle T}(t)\,B)\,R^{-1}\,(S^\tra\!+B^\tra P_{\scriptscriptstyle T}(t))+Q=0, \label{rde1} 
\end{eqnarray}
 with the terminal condition 
 \begin{equation}
 \label{tc}
 P_{\scriptscriptstyle T}(T)=0.  
 \end{equation}
If this differential equation admits solution $P_{\scriptscriptstyle T}(t)$ in $[0,T]$, then by following the same steps of
\cite[Theorem 3.1]{Ferrante-N-14}, we immediately see that
\begin{eqnarray*}
J_{T}(x_{\scriptscriptstyle 0},u) 
\ns&\ns = \ns&\ns \int_{\scriptscriptstyle 0}^T \| R^{-\frac{1}{2}}  (S^\tra+B^\tra P_{\scriptscriptstyle T}(t)) \,x(t) + R^{\frac{1}{2}}\,u(t) \|_2^2 \,dt +x^\tra(0)\,P_{\scriptscriptstyle T}(0)\,x(0),
\end{eqnarray*}
so that the optimal control is clearly
$u(t)=-R^{-1} (S^\tra+B^\tra P_{\scriptscriptstyle T}(t)) \,x(t)
$  and the optimal value of the cost is $J_{T}^\ast(x_{\scriptscriptstyle 0}) =x^\tra(0)\,P_{\scriptscriptstyle T}(0)\,x(0)$.
We now show that (\ref{rde1})-(\ref{tc}) indeed admit a unique solution  $P_{\scriptscriptstyle T}(t)$ in $(-\infty,T]$.
In fact, uniqueness is guaranteed by smoothness  of (\ref{rde1}) which also guarantees  existence of $P_{\scriptscriptstyle T}(t)$ in $(T-\varepsilon, T]$ for a sufficiently small $\varepsilon$.
To conclude it is therefore sufficient to show that no {\em finite escape time} can occur in this case. To this aim, consider
$P_{\scriptscriptstyle T}(T-t)=P_t(0)$ so that it is clear that as $t$ increases from zero to infinity,
$P_{\scriptscriptstyle T}(T-t)$ is bounded from below by the zero matrix, since $x^\tra(0) P_t(0)\,x(0)$ is the cost of a finite horizon LQ problem.
Moreover, since $R$ is positive definite, the solution $P_{\scriptscriptstyle T}(T-t)$ is also
bounded from above by the solution of the {final} value problem $\dot{P}_{ub}(t)=-[P_{\scriptscriptstyle T}(t)\,A+A^\tra\,P_{\scriptscriptstyle T}(t)+Q]$, $P_{ub}(T)=0$ in which there cannot be finite escape time because the differential equation is linear.
Thus, (\ref{rde1})-(\ref{tc}) admit a unique solution  $P_{\scriptscriptstyle T}(t)$ in $(-\infty,T]$.

Now consider the new matrix function 
$X(t)\defi P_t(0)=P_{\scriptscriptstyle T}(T-t)$, $t\ge 0$.
We immediately see that $X(t)$ satisfies equation (\ref{grde111}) with initial condition $X(0)=0$.
Moreover $X(t)$ is a non-decreasing flow of positive semidefinite matrices, i.e. $X(t+\delta t)\geq X(t)\geq 0$, for all $t,\delta t\geq 0$.
 We now show that $X(t)$ is a bounded function of $t \ge 0$. Indeed, given the $i$-th canonical basis vector $e_i$ of $\real^n$, {we have that for all $t\geq 0$, $e_i^\tra X(t) \,e_i=J_{t}^\ast(e_i)\le J_{\infty}(e_i,\overline{u}_i)$, where $\overline{u}_i$ is a control that renders $J_{\infty}(e_i,\overline{u}_i)$ finite, which exists by assumption.
 Therefore, $X(t)$ is non-decreasing and bounded, so that the limit $\overline{X}\defi\lim_{t\rightarrow\infty} X(t)$ exists and is finite}. Taking the limit on both sides of (\ref{grde111}) we immediately see that $\overline{X} \ge 0$ is indeed a solution of CARE($\Sigma$).
 Indeed, by repeating verbatim the same steps of \cite[Theorem 3.2]{Ferrante-N-14}, we see that
 $\overline{X}$ is the minimal positive semidefinite solution of CARE($\Sigma$) and that
 $u^\ast(t)=-R^{-1}(S^\top+B^\top\,\overline{X})\,x(t)$ is the optimal control.
 \endproof
}

As already observed, Proposition \ref{theprev} shows that the existence of symmetric positive semidefinite solutions of CGCARE($\Sigma$) guarantees that the associated LQ optimal control problem admits an impulse-free solution.

In order to claim that the solvability of CGCARE($\Sigma$) is {equivalent} to the fact that the LQ problem is solvable with non-impulsive control laws, the converse implication also needs to be proved. This is the task addressed in the following {result}, which proves the implication {\bf {\em (A)}} $\Rightarrow$ {\bf {\em (B)}} of Theorem \ref{main}.
%
\ \\
\begin{proposition}
Let the LQ problem admit a non-impulsive solution for every initial condition $x_{\scriptscriptstyle 0}\in \real^n$. Then, CGCARE($\Sigma$) admits a symmetric positive semidefinite solution.
\end{proposition}

\proof
Let the (possibly singular) LQ problem admit a non-impulsive solution for every initial condition $x_{\scriptscriptstyle 0}\in \real^n$. In view of \cite[Theorem 2]{Willems-KS-86}, the optimal control $u^\ast$ can be written as the static state feedback 
\bea
\label{app}
u^\ast(t)=-K\,x(t).
\eea
{
This result was given in  \cite{Willems-KS-86} under the assumption of stabilizability of the pair $(A,B)$. On the other hand, this assumption was only introduced to the end of exploiting \cite[Proposition 10]{Willems-KS-86}, dealing with the regular case, as taken from \cite[Theorem 6.1]{Hautus-S-83}. Lemma \ref{lemW} above generalizes \cite[Proposition 10]{Willems-KS-86} by just requiring the weaker assumption that the performance index $J_{\infty}(x_{\scriptscriptstyle 0},u)$ can be rendered finite from any initial condition $x_{\scriptscriptstyle 0}$ with a suitable control function $u(t)$, in place of the stabilizability of the pair $(A,B)$. Therefore, the proof of \cite[Theorem 2]{Willems-KS-86} can be carried out verbatim with just the assumption of the existence of a control that renders  $J_{\infty}(x_{\scriptscriptstyle 0},u)$ finite for any $x_{\scriptscriptstyle 0}\in \real^n$. }
By factorizing the Popov matrix as
\[
\bmat{cc} Q & S \\ S^\top & R \emat=\bmat{cc} C^\top \\ D^\top \emat \bmat{cc} C & D \emat,
\]
where $[\,C\;\;D\,]$ is of full row-rank, we can re-write (\ref{costinf}) as
 \begin{eqnarray}
 \label{costinf1}
J_\infty(x_{\scriptscriptstyle 0},u)&=&\int_{\scriptscriptstyle 0}^\infty y^\tra(t) y(t)\,dt,
\end{eqnarray}
where $y(t)=C\,x(t)+D\,u(t)$ can be considered as a fictitious output function. 
{The closed-loop} system that corresponds to the application of the control (\ref{app}) is
\beann
\left\{ \begin{array}{lll} 
\dot{x}(t)\ns&\ns=\ns&\ns (A-B\,K)\,x(t) \\
y(t)\ns&\ns=\ns&\ns (C-D\,K)\,x(t) \end{array} \right.
\eeann
Let $A_{\scriptscriptstyle K} \defi A-B\,K$ and $C_{\scriptscriptstyle K} \defi C-D\,K$. The optimal state is
$x(t)=e^{A_{\scriptscriptstyle K}\,t}\,x_{\scriptscriptstyle 0}$,
and the corresponding output is 
$y(t)=C_{\scriptscriptstyle K}\,e^{A_{\scriptscriptstyle K}\,t}\,x_{\scriptscriptstyle 0}$. Thus, the optimal
cost is given by
\[
J_\infty(x_{\scriptscriptstyle 0},u^\ast)=x_{\scriptscriptstyle 0}^\top\,\left[ \int_{\scriptscriptstyle 0}^{\infty} e^{A_{\scriptscriptstyle K}^\top\,t}\,C_{\scriptscriptstyle K}^\top\,C_{\scriptscriptstyle K}\,e^{A_{\scriptscriptstyle K}\,t}\,dt\right]\,x_{\scriptscriptstyle 0}
\]
Let $r$ be the rank of $R$. Consider a basis of the input space such that
\[
D=[\,D_1\;\;0\,]\qquad \text{and} \qquad B=[\,B_1\;\;B_2\,],
\]
where $D_1$ is of full column-rank $r$. In this basis, we have $R=\bsmat R_1 & 0 \\[1mm] 0 & 0 \esmat$ and $S=\bsmat S_1 & 0 \esmat$, where $R_1 \in \real^{r \times r}$ is invertible and $S_1$ has $r$ columns. Let us now consider $x_{\scriptscriptstyle 0} \in \ima B_2$. Using a control $u^\circ= \bsmat 0_r \\[1mm] u^\circ_2 \esmat$ such that $u^\circ_2(t)$ is allowed to contain impulses (i.e., Dirac deltas and its derivatives in the distributional sense), the state can be instantaneously driven to the origin, i.e., $x(0^+)=0$, and $J_\infty(x_{\scriptscriptstyle 0},u^\ast)=0$ because in this basis the second block of components of the control law are not weighted in the performance index. Thus, 
$\ima B_2 \subseteq  \ker (C_{\scriptscriptstyle K}\,e^{A_{\scriptscriptstyle K}\,t})$,
so that 
\bea
\label{arb}
C_{\scriptscriptstyle K}\,e^{A_{\scriptscriptstyle K}\,t}\,B_2=0 \quad \forall\,t \ge 0,
\eea
which means that the transfer function $C_{\scriptscriptstyle K}\,(s\,I_n-A_{\scriptscriptstyle K})^{-1}\,B_2$ is zero.
 Let $x_{\scriptscriptstyle 0} \in \real^n$, and $u^\ast$ be a corresponding optimal control. Let $u^\ast$ be partitioned as 
$u^\ast(t)=\bsmat u_1^\ast(t) \\[1mm] u_2^\ast(t) \esmat$, 
conformably with the decomposition of the input space. Then, given any $\delta\,u_2(t)$, we can define the new input
$\tilde{u}^\ast(t) \defi \bsmat u_1^\ast(t) \\[1mm] u_2^\ast(t)+\delta\,u_2(t) \esmat$. 
Thus, (\ref{arb}) guarantees that 
$y_{x_{\scriptscriptstyle 0},u^\ast}(t) = y_{x_{\scriptscriptstyle 0},\tilde{u}^\ast}(t)$, where $y_{x_{\scriptscriptstyle 0},u^\ast}(t)$ is the output that corresponds to $x_{\scriptscriptstyle 0}$ and $u^\ast$ while $y_{x_{\scriptscriptstyle 0},\tilde{u}^\ast}(t)$ is the one that corresponds to $x_{\scriptscriptstyle 0}$ and $\tilde{u}^\ast$, this in turn implies that
$J^\star \defi J(x_{\scriptscriptstyle 0},u^\ast)=J(x_{\scriptscriptstyle 0},\tilde{u}^\ast)$. 
Hence, the (regular) LQ problem for the quadruple $(A,B_1,C,D_1)$, i.e., the one consisting of the minimization of the performance index
\[
\hat{J}(x_{\scriptscriptstyle 0},u_1) \defi \int_{\scriptscriptstyle 0}^\infty  \bmat{cc} x^\tra(t) & u_1^\tra(t) \emat \bmat{cc} Q & S_1 \\ S_1^\tra & R_1 \emat \bmat{c} x(t) \\ u_1(t) \emat\,dt
\]
subject to the constraint 
 $\dot{x}(t)=A\,x(t)+B_1\,u_1(t)$ and $x(0)=x_{\scriptscriptstyle 0}$, admits solutions for all $x_{\scriptscriptstyle 0}$, and the corresponding optimal cost coincides with the optimal cost of the original LQ problem, which is $\hat{J}(x_{\scriptscriptstyle 0},u_1^\ast)=J^\star$. On the other hand, as already observed, since $R_1=D_1^\top D_1$ is positive definite, this LQ problem for the quadruple $(A,B_1,C,D_1)$ is regular. The fact that it admits solutions for all $x_{\scriptscriptstyle 0}$ implies that the corresponding algebraic Riccati equation 
\bea
\label{CARE1}
X\,A+A^\top\,X-(C^\top D_1+X\,B_1)\,(D_1^\top D_1)^{-1} (D_1^\top C+B_1^\top X)+C^\top C=0
\eea
 admits a solution $\overline{X}=\overline{X}^\top \ge 0$, and $J^\star=x_{\scriptscriptstyle 0}^\top\,\overline{X}\,x_{\scriptscriptstyle 0}$. Thus,
\bea
\label{new}
\overline{X}=\int_{\scriptscriptstyle 0}^\infty e^{A_{\scriptscriptstyle K}^\top\,t}\,C_{\scriptscriptstyle K}^\top C_{\scriptscriptstyle K} \,e^{A_{\scriptscriptstyle K}\,t}\,dt.
\eea
We can re-write (\ref{CARE1}) in the form
\[
X\,A+A^\top\,X-\bmat{cc} C^\top D_1+X\,B_1 & X\,B_2 \emat
\bmat{cc} D_1^\top D_1 & 0 \\ 0 & 0 \emat \bmat{c} 
D_1^\top C+B_1^\top X \\ B_2^\top X \emat+C^\top C =0,
\]
which is exactly the original GCARE($\Sigma$)
\[
X\,A+A^\top\,X -(C^\top D+X\,B)\,(D^\top D)^{\dagger} (D^\top C+B^\top X)+C^\top C=0
\]
Thus, $\overline{X}=\overline{X}^\top\ge0$ is a solution of GCARE($\Sigma$). Moreover, from (\ref{arb}) we have $\ima B_2 \subseteq \ker (C_{\scriptscriptstyle K}\,e^{A_{\scriptscriptstyle K}\,t})$ for all $t \ge 0$, which, together with (\ref{new}), yields
$\ima B_2 \subseteq \ker \overline{X}$.
It is easy to see that this means that 
$\ker R \subseteq \ker (S+\overline{X}\,B)$.
Indeed, in the chosen basis this subspace inclusion reads as
\[
\ima \bmat{c} 0 \\ I \emat =\ker \bmat{cc} D_1^\top\,D_1 & 0 \\ 0 & 0 \emat \subseteq \ker \bmat{cc} C\,D_1+\overline{X}\,B_1 & X\,B_2 \emat = \ker \bmat{cc} C\,D_1+\overline{X}\,B_1 & 0 \emat,
\]
which is certainly satisfied. Thus, $\overline{X}$ is also a symmetric and positive semidefinite solution of CGCARE($\Sigma$).
\endproof

{
Notice that, as a byproduct  of the previous proof, in the so-called {\em cheap} case, i.e. when $R=0$, we have the following 
\begin{corollary} 
Let $R=0$. If Problem \ref{prb1} admits a regular solution for any initial condition $x_{\scriptscriptstyle 0}$ then
the optimal cost is zero: $J^\star(x_{\scriptscriptstyle 0})=0$ for each $x_{\scriptscriptstyle 0}\in\real^n$. 
\end{corollary}
}

\subsection{Geometric conditions}
So far, we have proved that the statements  {\bf {\em (A)}}, {\bf {\em (B)}} and {\bf {\em (C)}} in Theorem \ref{main} are equivalent. In this section, we focus our attention on condition {\bf {\em (D)}} of the same theorem, and we show that it is also equivalent to the other three conditions.

Consider the quadruple $(A,B,C,D)$, where $C$ and $D$ are matrices of suitable sizes such that
(\ref{fact}) holds.

{\begin{proposition}
\label{progeom}
Let CGCARE($\Sigma$) admit a solution $X=X^\tra$. Then, $\gS^\star=\gR^\star$.
\end{proposition}
}
\proof
Let  $X=X^\tra$ be a solution of CGCARE($\Sigma$). Observe also that CGCARE($\Sigma$) can be re-written as
\begin{equation}
\label{cgdare12}
\left\{ \begin{array}{ll} X\,A_{\scriptscriptstyle 0}+A_{\scriptscriptstyle 0}^\tra\,X-X\,B\,R^\dagger\,B^\tra\,X+Q_{\scriptscriptstyle 0}=0 \\
\ker R \subseteq \ker X\,B \end{array} \right.
\end{equation}
where $A_{\scriptscriptstyle 0} \defi A-B\,R^\dagger S^\tra$ and $Q_{\scriptscriptstyle 0} \defi Q-S\,R^\dagger S^\tra$. Recall that $G=I_m-R^\dagger R$, so that $B_2 = B\,G$, and (\ref{cgdare12}) becomes
\begin{equation}
\label{cgdare13}
\left\{ \begin{array}{ll} X\,A_{\scriptscriptstyle 0}+A_{\scriptscriptstyle 0}^\tra\,X-X\,B\,R^\dagger\,B^\tra\,X+Q_{\scriptscriptstyle 0}=0 \\
X\,B\,G=0 \end{array} \right.
\end{equation}
It is easy to see that $\ker X \subseteq \ker Q_{\scriptscriptstyle 0}$. 
Indeed, by multiplying the first of (\ref{cgdare13}) on the left by $\xi^\tra$ and on the right by $\xi$, where $\xi \in \ker X$, we get $\xi^\tra\,Q_{\scriptscriptstyle 0}\,\xi=0$. However, $Q_{\scriptscriptstyle 0}$ is positive semidefinite, being the generalized Schur complement of $Q$ in $\Pi$. Hence, $Q_{\scriptscriptstyle 0}\,\xi=0$, which implies $\ker X \subseteq \ker Q_{\scriptscriptstyle 0}$. 
Since $X\,B\,G=0$, we get also $Q_{\scriptscriptstyle 0}\,B\,G=0$. By post-multiplying the first of (\ref{cgdare13}) by a vector $\xi \in \ker X$ we find $X\,A_{\scriptscriptstyle 0}\,\xi=0$, which says that $\ker X$ is $A_{\scriptscriptstyle 0}$-invariant. This means that $\ker X$ is an $A_{\scriptscriptstyle 0}$-invariant subspace containing the image of $B\,G$. Then, the reachable subspace of the pair $(A_{\scriptscriptstyle 0},BG)$, denoted by $\gR(A_{\scriptscriptstyle 0},B\,G)$, which is the smallest $A_{\scriptscriptstyle 0}$-invariant subspace containing the image of $B\,G$, is contained in $\ker X$, i.e., $\gR(A_{\scriptscriptstyle 0},B\,G)\subseteq \ker X$. Therefore also $\gR(A_{\scriptscriptstyle 0},B\,G)\subseteq \ker Q_{\scriptscriptstyle 0}$. Notice that $Q_{\scriptscriptstyle 0}$ can be written as $C_{\scriptscriptstyle 0}^\tra\,C_{\scriptscriptstyle 0}$, where $C_{\scriptscriptstyle 0}\defi C-D\,R^\dagger S^\tra$. Indeed,
\begin{eqnarray*}
C_{\scriptscriptstyle 0}^\tra\,C_{\scriptscriptstyle 0} &=& C^\tra C-C^\tra D R^\dagger S^\tra-S R^\dagger D^\tra\,C+S R^\dagger D^\tra D R^\dagger S^\tra \\ 
&= &Q-S R^\dagger S-S R^\dagger S^\tra+S R^\dagger S^\tra=Q_{\scriptscriptstyle 0}.
\end{eqnarray*}
Consider the two quadruples $(A,B,C,D)$ and $(A_{\scriptscriptstyle 0},B,C_{\scriptscriptstyle 0},D)$. We observe that the second is obtained directly from the first by applying the feedback input $u(t)=-R^\dagger S\,x(t)+v(t)$. We denote by $\gV^\star$, $\gR^\star$ the largest output-nulling and reachability subspace of $(A,B,C,D)$, and by $\gS^\star$ the smallest input-containing subspace of $(A,B,C,D)$. Likewise, we denote by $\gV_{\scriptscriptstyle 0}^\star$, $\gR_{\scriptscriptstyle 0}^\star$, $\gS_{\scriptscriptstyle 0}^\star$ the same subspaces relative to the quadruple $(A_{\scriptscriptstyle 0},B,C_{\scriptscriptstyle 0},D)$. Thus, $\gV^\star= \gV^\star_{\scriptscriptstyle 0}$, 
$\gR^\star= \gR^\star_{\scriptscriptstyle 0}$, and $\gS^\star= \gS^\star_{\scriptscriptstyle 0}$. 
 The first two identities are obvious, since output-nulling subspaces can be made invariant under state-feedback transformations and reachability is invariant under the same transformation. The third follows from \cite[Theorem 8.17]{Trentelman-SH-01}.
 There holds $\gR^\star=\gR(A_{\scriptscriptstyle 0},B\,G)$. Indeed, consider a state $x_1\in \gR(A_{\scriptscriptstyle 0},B\,G)$. There exists a control function $u$ driving the state from the origin to $x_1$, and we show that this control keeps the output at zero. Since $\ima(B\,G)=B\,\ker D$, such control can be chosen to satisfy $D\,u(t)=0$ for all $t\ge 0$. Moreover, as we have already seen, from $Q_{\scriptscriptstyle 0}=C_{\scriptscriptstyle 0}^\tra C_{\scriptscriptstyle 0}$ and $\gR(A,BG)=\gR(A_{\scriptscriptstyle 0},BG)$ we have $C_{\scriptscriptstyle 0}\,\gR(A_{\scriptscriptstyle 0},B\,G)=0$ since $\gR(A,B\,G)$ lies in $\ker Q_{\scriptscriptstyle 0}$. Therefore, the output is identically zero. This implies that $\gR(A_{\scriptscriptstyle 0},B\,G) \subseteq \gR^\star$. However, the reachability subspace of $(A_{\scriptscriptstyle 0},B,C_{\scriptscriptstyle 0},D)$ cannot be greater than $\gR(A_{\scriptscriptstyle 0},B\,G)$, since 
 $D^\tra C_{\scriptscriptstyle 0}=D^\tra (I_m-D\,(D^\tra D)^{\dagger}D^\tra)C=0$. Therefore, such control must necessarily render the output non-zero. The same argument can be used to prove that $\gS^\star=\gR(A_{\scriptscriptstyle 0},B\,G)$, where distributions can also be used in the allowed control, since $\gR(A,BG)$ represents also the set of states that are reachable from the origin using distributions in the control law \cite[p. 183]{Trentelman-SH-01}. Hence, $\gS^\star=\gR^\star$.
 \endproof

\begin{remark}
Proposition \ref{progeom} proves a stronger result than the implication of {\bf {\em (C)}} $\Rightarrow$ {\bf {\em (D)}} in Theorem \ref{main}. On the other hand, it is easy to see that the converse of this result does not hold, unless we introduce -- as in Theorem \ref{main} -- the additional assumption that for every initial state the performance index can be made finite. Indeed, consider an LQ problem where
\[
A=\bmat{cc} 0 & 0 \\ 0 & 1 \emat, \quad B=\bmat{c} 0 \\ 1 \emat, \quad Q=\bmat{cc} 1 & 0 \\ 0 & 0 \emat, 
\]
and $S=0$ and $R=0$, so that $C=[\,1\;\;\;0\,]$ and $D=0$. In this case, it is found that 
\[
\gV^\star=\gS^\star=\gR^\star=\spanR \left\{ \bmat{c} 0 \\ 1 \emat \right\},
\]
In this case the CGCARE($\Sigma$) reduces to the Lyapunov equation $X\,A+A^\top\,X+Q=0$. Partitioning $X$ as
$X=\bsmat x_1 & x_2 \\[1mm] x_2 & x_3 \esmat$, 
the Lyapunov equation becomes
\[
\bmat{cc} 1 & x_2 \\ x_2 & 2\,x_3 \emat=0,
\]
which clearly does not admit solutions. However, it is easily seen that in this example the state dynamics are
\beann
\dot{x}_1(t) \ns&\ns = \ns&\ns  0 \\
\dot{x}_2(t) \ns&\ns = \ns&\ns  x_2(t)+u(t) 
\eeann
and the performance index is $J_{\infty}(x_{\scriptscriptstyle 0},u)=\int_{\scriptscriptstyle 0}^\infty x_1^2(t) \,dt$, which is not finite if $x_1(0) \neq 0$.
\end{remark}

{The following result shows that {\bf {\em (D)}} $\Rightarrow$ {\bf {\em (A)}}, completing the proof} of Theorem \ref{main}.

\begin{proposition}
\label{progeom-conv}
Let $\gS^\star=\gR^\star$, and assume that for every initial condition $x_{\scriptscriptstyle 0}$ there exists a control $u$ such that $J_{\infty}(x_{\scriptscriptstyle 0},u)$ is finite. Then, 
there exists a non-impulsive optimal control.
\end{proposition}
\proof
Let $\gS^\star=\gR^\star$. 
Consider the decomposition in \cite[p. 328]{Willems-KS-86}. If $\gS^\star=\gR^\star$, the fourth and the fifth block components of the state disappear, and the system dynamics reduce to
\beann
\bmat{c}
\dot{x}_1(t) \\
\dot{x}_2(t) \\ 
\dot{x}_3(t) \emat \ns&\ns = \ns&\ns \bmat{ccc} A_{11} & 0 & 0 \\ A_{21} & A_{22} & 0 \\
0 & A_{32} & A_{33} \emat \bmat{c} x_1(t) \\x_2(t) \\x_3(t) \emat +
\bmat{c} B_{11} \\ B_{12} \\ B_{13} \emat u_1^\prime(t)+ 
\bmat{c} 0 \\ 0 \\ B_{13} \emat u_2^\prime(t) \\
y_1(t) \ns&\ns = \ns&\ns u_1^\prime(t) \\
y_2(t) \ns&\ns = \ns&\ns \bmat{ccccc} C_{21} & 0 & 0 \emat \bmat{c} x_1(t) \\x_2(t) \\x_3(t) \emat
\eeann
In view of \cite[Theorem 2]{Willems-KS-86}, the only part of the state where there may be distributions in the optimal control is the third. On the other hand, the third block of coordinates of this basis span $\gR^\star$. This implies that $x_3$ is arbitrary, in the sense that it is not penalized in the performance index. Thus, an optimal control such that there are distributions in $x_3$ continues to be optimal even when such distributions are removed. Therefore, the optimal control can be rendered regular.
\endproof

 \section{Concluding remarks}
    In this paper, a full picture has been drawn illustrating the relationship that exists between the solvability of the so-called constrained generalized Riccati equation and the existence of non-impulsive optimal controls of the associated infinite-horizon LQ problem. This link has been examined both from an algebraic and a geometric angle. Now that this relationship has been clarified and explained, an important direction of future research aims at obtaining a full characterization of the set of solutions of the constrained generalized continuous algebraic Riccati 
equation that parallels the discrete time counterpart in \cite{Ferrante-N-12-sub,Ferrante-N-13-1}.

\end{document}